%% file: ejcarticle.tex
\documentclass[12pt]{article}

\makeatletter
\newfont{\footsc}{cmcsc10 at 8truept}
\newfont{\footbf}{cmbx10 at 8truept}
\newfont{\footrm}{cmr10 at 10truept}
\renewcommand{\ps@plain}{%
\renewcommand{\@oddfoot}{\footsc the electronic journal of combinatorics
  {\footbf 12} (2005), \#R00\hfil\footrm\thepage}}
\makeatother 

\input{headerejc.tex}
\usepackage{amssymb}
\usepackage{young}
\usepackage[vcentermath]{youngtab}
\begin{document}

\begin{center} \Large \textbf{The hook fusion procedure}\\
\quad \\
\normalsize James Grime \\ \small \emph{Department of Mathematics,
University of York, York, YO10 5DD, UK} \\
\emph{jrg112@york.ac.uk}
\end{center}

\small Submitted: Feb 25, 2005;  Accepted: Apr 25, 2005; Published: June 25, 2005

\quad

\textbf{Abstract}\\ We derive a new expression for the diagonal
matrix elements of irreducible representations of the symmetric
group. We obtain this new expression using Cherednik's fusion
procedure. However, instead of splitting Young diagrams into their
rows or columns, we consider their principal hooks. This minimises
the number of auxiliary parameters needed in the fusion procedure.

\section{Introduction}

In this article we present a new expression for the diagonal
matrix elements of irreducible representations of the symmetric
group. We will obtain this new expression using Cherednik's fusion
procedure. This method originates from the work of Jucys
\cite{5a}, and has already been used by Nazarov and Tarasov
\cite{7a, 9}. However our approach differs by minimizing the
number of auxiliary parameters needed in the fusion procedure.
This is done by considering hooks of Young diagrams, rather than
their rows or columns as in \cite{7a, 9}.

Irreducible representations of the group $S_n$ are parameterized
by partitions of $n$, \cite{3}. The \emph{Young diagram} of a
partition $\lambda$ is the set of boxes $(i, j) \in \mathbb{Z}^2$
such that $1 \leqslant j \leqslant \lambda_i$. In drawing such
diagrams we let the first coordinate $i$ increase as one goes
downwards, and the second coordinate $j$ increase from left to
right. For example the partition $\lambda = (3,3,2)$ gives the
diagram \[ \yng(3,3,2)\]

If $(i,j)$ is a box in the diagram of $\lambda$, then the
$(i,j)$-\emph{hook} is the set of boxes in $\lambda$
\[  \{ (i, j') : j' \geqslant j \} \cup \{ (i', j)
: i' \geqslant i \}, \] We call the $(i,i)$-hook the
\emph{$i^{\mbox{th}}$ principal hook}.

A \emph{standard tableau}, $\Lambda$, is a filling of the diagram
$\lambda$ in which the entries are the numbers 1 to $n$, each
occurring once. The Young Symmetrizer of a standard tableau
$\Lambda$ of shape $\lambda$ generates an irreducible
$\mathbb{C}S_n$-module, \cite{3}, denoted $V_\lambda$.

In 1986 Ivan Cherednik proposed another description of the Young
symmetrizer, \cite{2}. The following version of Cherednik's
description can be found in \cite{9}.

There is a basis for the space $V_\lambda$ called the \emph{Young
basis}, with its vectors labeled by the standard tableaux of shape
$\lambda$, \cite{6a}, \cite{10}. Now, for any standard tableau
$\Lambda$ consider the \emph{diagonal matrix element},
$F_\Lambda$, of the representation $V_\lambda$ corresponding to
the Young basis vector $v_\Lambda$, defined by,
\[ F_\Lambda = \sum_{g \in S_n} (v_\Lambda, g v_\Lambda)
g \in \mathbb{C}S_n. \]

Multiplying $F_\Lambda$ on the left and right by certain
invertible elements returns the Young symmetrizer. Hence
$F_\Lambda$ may also be used to generate the irreducible module
$V_\lambda$, \cite{6a}, \cite[section 2]{8}.

Cherednik's fusion procedure can be used to find an alternative
expression for $F_\Lambda$. Consider a standard tableaux $\Lambda$
of shape $\lambda$. For each $p = 1, \dots, n$, put $c_p = j-i$ if
the box $(i,j) \in \lambda$ is filled with the number $p$ in
$\Lambda$. The difference $j-i$ is the \emph{content} of box
$(i,j)$.

For any two distinct numbers $p$, $q \in \{1, \dots, n\}$, let
$(pq)$ be the transposition in the symmetric group $S_n$. Consider
the rational function of two complex variables $u$, $v$, with
values in the group ring $\mathbb{C}S_n$:
\begin{equation}\label{smallf} f_{pq}(u,v) = 1 -
\frac{(pq)}{u-v}.\end{equation}

Now introduce $n$ complex variables $z_1, \dots , z_n$. Consider
the set of pairs $(p,q)$ with \\$1 \leqslant p < q \leqslant n$.
Ordering the pairs lexicographically we form the product
\begin{equation}\label{bigf} F_\Lambda(z_1, \dots, z_n) = \prod_{(p,q)}^\rightarrow
f_{pq}(z_p + c_p, z_q + c_q).\end{equation} If $p$ and $q$ sit on
the same diagonal in the tableau $\Lambda$, then $f_{pq}(z_p +
c_p, z_q + c_q)$ has a pole at $z_p = z_q$.

Let $\mathcal{R}_\Lambda$ be the vector subspace in $\mathbb{C}^n$
consisting of all tuples $(z_1, \dots , z_n)$ such that $z_p =
z_q$ whenever the numbers $p$ and $q$ appear in the same row of
the tableau $\Lambda$.

By direct calculation we can show the following the identities are
true; \begin{equation}\label{triple}
f_{pq}(u,v)f_{pr}(u,w)f_{qr}(v,w) =
f_{qr}(v,w)f_{pr}(u,w)f_{pq}(u,v) \end{equation} for all pairwise
distinct indices $p$, $q$, $r$, and
\begin{equation}\label{commute} f_{pq}(u,v)f_{st}(z,w) =
f_{st}(z,w)f_{pq}(u,v) \end{equation} for all pairwise distinct
$p$, $q$, $s$, $t$.

We call (\ref{triple}) and (\ref{commute}) the Yang-Baxter
relations. Using these relations we obtain reduced expressions for
$F_\Lambda(z_1, \dots, z_n)$ different from the right hand side of
(\ref{bigf}). For more details on different expressions for
$F_\Lambda(z_1, \dots, z_n)$ see \cite{5aa}.

Using (\ref{triple}) and (\ref{commute}) we may reorder the
product $F_\Lambda(z_1, \dots , z_n)$ such that each singularity
is contained in an expression known to be regular at $z_1 = z_2 =
\dots = z_n$, \cite{7}. It is by this method that it was shown
that the restriction of the rational function $F_\Lambda(z_1,
\dots ,z_n)$ to the subspace $\mathcal{R}_\Lambda$ is regular at
$z_1 = z_2 = \dots = z_n$. Furthermore, by showing divisibility on
the left and on the right by certain elements of $\mathbb{C}S_n$
it was shown that the value of $F_\Lambda(z_1, \dots ,z_n)$ at
$z_1 = z_1 = \dots = z_n$ is the diagonal matrix element
$F_\Lambda$, \cite{3}. Thus we have the following theorem,

\begin{nottheorem} \label{MNtheorem} Restriction to $\mathcal{R}_\Lambda$ of the
rational function $F_\Lambda(z_1, \dots , z_n)$ is regular at $z_1
= z_2 = \dots = z_n$. The value of this restriction at $z_1 = z_2
= \dots = z_n$ coincides with the element $F_\Lambda \in
\mathbb{C}S_n$.
\end{nottheorem}

Similarly, we may form another expression for $F_\Lambda$ by
considering the subspace in $\mathbb{C}^n$ consisting of all
tuples $(z_1, \dots , z_n)$ such that $z_p = z_q$ whenever the
numbers $p$ and $q$ appear in the same column of the tableau
$\Lambda$ \cite{9}.

In this article we present a new expression for the diagonal
matrix elements which minimizes the number of auxiliary parameters
needed in the fusion procedure. We do this by considering hooks of
standard tableaux rather than their rows or columns.

Let $\mathcal{H}_\Lambda$ be the vector subspace in $\mathbb{C}^n$
consisting of all tuples $(z_1, \dots , z_n)$ such that $z_p =
z_q$ whenever the numbers $p$ and $q$ appear in the same principal
hook of the tableau $\Lambda$. We will prove the following
theorem.

\begin{theorem}\label{fulltheorem} Restriction to $\mathcal{H}_\Lambda$ of the
rational function $F_\Lambda(z_1, \dots , z_n)$ is regular at $z_1
= z_2 = \dots = z_n$. The value of this restriction at $z_1 = z_2
=  \dots = z_n$ coincides with the element $F_\Lambda \in
\mathbb{C}S_n$.
\end{theorem}

In particular, this hook fusion procedure can be used to form
irreducible representations of $S_n$ corresponding to Young
diagrams of hook shape using only one auxiliary parameter, $z$. By
taking this parameter to be zero we find that no parameters are
needed for diagrams of hook shape. Therefore if $\nu$ is a
partition of hook shape we have
\begin{equation}\label{bigfhook} F_\nu = F_\nu(z) =
\prod_{(p,q)}^\rightarrow f_{pq} (c_p, c_q), \end{equation} with
the pairs $(p, q)$ in the product ordered lexicographically.


To motivate the study of modules corresponding to partitions of
hook shape first let us consider  the \emph{Jacobi-Trudi
identities} \cite[Chapter I3]{6}. We can think of the following
identities as dual to the Jacobi-Trudi identities.

If $\lambda = (\lambda_1, \dots , \lambda_k)$ such that $n =
\lambda_1 + \dots + \lambda_k$ then we have the following
decomposition of the induced representation of the tensor product
of modules corresponding to the rows of $\lambda$;
\[ \textrm{Ind}_{S_{\lambda_1} \times S_{\lambda_2} \times \cdots
\times S_{\lambda_k}}^{S_n } V_{(\lambda_1)} \otimes
V_{(\lambda_2)} \otimes \cdots \otimes V_{(\lambda_k)} \cong
\bigoplus_{\mu} (V_\mu)^{\oplus K_{\mu \lambda}},\] where the sum
is over all partitions of $n$. Note that $V_{(\lambda_i)}$ is the
trivial representation of $S_{\lambda_i}$. The coefficients
$K_{\mu \lambda}$ are non-negative integers known as \emph{Kostka
numbers}, \cite{6}. Importantly, we have $K_{\lambda \lambda} =
1.$

On the subspace $\mathcal{R}_\Lambda$, if $z_i - z_j \notin
\mathbb{Z}$ when $i$ and $j$ are in different rows of $\Lambda$
then the above induced module may be realised as the left ideal in
$\mathbb{C}S_n$ generated by $F_\Lambda(z_1, \dots, z_n)$.\\
The irreducible representation $V_\lambda$ appears in the
decomposition of this induced module with coefficient 1, and is
the ideal of $\mathbb{C}S_n$ generated by $F_\Lambda(z_1, \dots ,
z_n)$ when $z_1 = z_2 = \dots = z_n$. The fusion procedure of
theorem \ref{MNtheorem} provides a way of singling out this
irreducible component.

Similarly we have the equivalent identity for columns, \[
\textrm{Ind}_{S_{\lambda'_1} \times S_{\lambda'_2} \times \cdots
\times S_{\lambda'_l}}^{S_n } V_{(1^{\lambda'_1})} \otimes
V_{(1^{\lambda'_2})} \otimes \cdots \otimes V_{(1^{\lambda'_l})}
\cong \bigoplus_{\mu} (V_{\mu})^{\oplus K_{\mu' \lambda'}},\]
where $l$ is the number of columns of $\lambda$. In this case
$V_{(1^{\lambda'_i})}$ is the alternating representation of
$S_{\lambda'_i}$. This induced module is isomorphic to the left
ideal of $\mathbb{C}S_n$ generated by $F_{\Lambda}(z_1, \dots,
z_n)$ considered on the subspace $\mathcal{R}_{\Lambda'}$, with
$z_i - z_j \notin \mathbb{Z}$ when $i, j$ are in different columns
of $\Lambda$. Again the irreducible representation $V_{\lambda}$
appears in the decomposition of this induced module with
coefficient 1, and is the ideal of $\mathbb{C}S_n$ generated by
$F_{\Lambda}(z_1, \dots, z_n)$ when $z_1 = z_2 = \dots = z_n$.

There is another expression known as the \emph{Giambelli identity}
\cite{5}. Unlike the Jacobi-Trudi identities, this identity
involves splitting $\lambda$ into its principal hooks, rather than
its rows or columns. A combinatorial proof of the Giambelli
identity can be found in \cite{4}.

Divide a Young diagram $\lambda$ into boxes with positive and
non-positive content. We may illustrate this on the Young diagram
by drawing 'steps' above the main diagonal. Denote the boxes above
the steps by $\alpha(\lambda)$ and the rest by $\beta(\lambda)$.
For example, the following figure illustrates $\lambda$,
$\alpha(\lambda)$ and $\beta(\lambda)$ for $\lambda = (3,3,2)$.

\begin{figure}[h]
\label{steps}
\begin{center}
\begin{picture}(150, 60)
\put(0,30){\framebox(15,15)[r]{  }}
\put(15,30){\framebox(15,15)[r]{ }}
\put(30,30){\framebox(15,15)[r]{  }}

\put(0,15){\framebox(15,15)[r]{  }}
\put(15,15){\framebox(15,15)[r]{  }}
\put(30,15){\framebox(15,15)[r]{  }}

\put(0,0){\framebox(15,15)[r]{  }}
\put(15,0){\framebox(15,15)[r]{}}

\put(72.5,30){\framebox(15,15)[r]{  }}
\put(87.5,30){\framebox(15,15)[r]{ }}
\put(125,30){\framebox(15,15)[r]{ }}

\put(87.5,15){\framebox(15,15)[r]{  }}
\put(125,15){\framebox(15,15)[r]{  }}
\put(140,15){\framebox(15,15)[r]{  }}

\put(125,0){\framebox(15,15)[r]{ }}
\put(140,0){\framebox(15,15)[r]{}}

\put(20,50){$\lambda$} \put(77.5,50){$\alpha(\lambda)$}
\put(130,50){$\beta(\lambda)$}

\linethickness{2pt} \put(0,45){\line(1,0){15}}
\put(15,30){\line(1,0){15}} \put(30,15){\line(1,0){15}}

\put(15,30){\line(0,1){15}} \put(30,15){\line(0,1){15}}

\end{picture}
\end{center}
 \caption{The Young diagram $(3,3,2)$ divided into boxes with positive content and non-positive content}
 \end{figure}
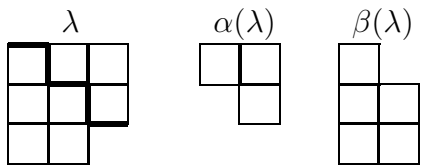

If we denote the rows of $\alpha(\lambda)$ by $\alpha_1 > \alpha_2
> \dots > \alpha_d > 0$ and the columns of $\beta(\lambda)$
by $\beta_1 > \beta_2 > \dots > \beta_d > 0$, then we have the
following alternative notation for $\lambda$;
\[ \lambda = ( \alpha | \beta ), \]
where $\alpha = (\alpha_1, \dots , \alpha_d)$ and $\beta =
(\beta_1, \dots , \beta_d)$.

Here $d$ denotes the length of the side of the \emph{Durfee
square} of shape $\lambda$, which is the set of boxes
corresponding to the largest square that fits inside $\lambda$,
and is equal to the number of principal hooks in $\lambda$. In our
example $d=2$ and $\lambda = (2, 1 | 3, 2)$.

We may consider the following identity as a dual of the Giambelli
identity.
\[ \textrm{Ind}_{S_{h_1} \times S_{h_2} \times \cdots
\times S_{h_d}}^{S_n} V_{(\alpha_1 | \beta_1)} \otimes
V_{(\alpha_2 | \beta_2)} \otimes \cdots \otimes V_{(\alpha_d |
\beta_d)} \cong \bigoplus_{\mu} (V_{\mu})^{\oplus D_{\mu
\lambda}},\] where $h_i$ is the length of the $i^{\mbox{th}}$
principal hook, and the sum is over all partitions of $n$. This is
a decomposition of the induced representation of the tensor
product of modules of hook shape. Further these hooks are the
principal hooks of $\lambda$. The coefficients, $D_{\mu \lambda}$,
are non-negative integers, and in particular $D_{\lambda \lambda}
=1$.

On the subspace $\mathcal{H}_\Lambda$, if $z_i - z_j \notin
\mathbb{Z}$ when $i$ and $j$ are in different principal hooks of
$\Lambda$ then the above induced module may be realised as the
left ideal in
$\mathbb{C}S_n$ generated by $F_\Lambda(z_1, \dots, z_n)$.\\
The irreducible representation $V_\lambda$ appears in the
decomposition of this induced module with coefficient 1, and is
the ideal of $\mathbb{C}S_n$ generated by $F_\Lambda(z_1, \dots ,
z_n)$ when $z_1 = z_2 = \dots = z_n$.

Hence, in this way,  our hook fusion procedure relates to the
Giambelli identity in the same way that Cherednik's original
fusion procedure relates to the Jacobi-Trudi identity. Namely, it
provides a way of singling out the irreducible component
$V_\lambda$ from the above induced module.

The fusion procedure was originally developed in the study of
affine Hecke algebras, \cite{2}. Our results may be regarded as an
application of the representation theory of these algebras,
\cite{10}.


Acknowledgements and thanks go to Maxim Nazarov for his
supervision, and for introducing me to this subject. I would also
like to thank EPSRC for funding my research and my anonymous
referees for their comments.

\section{Fusion Procedure for a Young Diagram}

The diagonal matrix element $F_\Lambda$ determines the irreducible
module $V_\lambda$ of $S_n$, up to isomorphism, for any tableau
$\Lambda$ of shape $\lambda$. Therefore in the sequel we will only
use one particular tableau, the \emph{hook tableau}. In which case
we may denote the diagonal matrix element $F_\Lambda$ by
$F_\lambda$, and the space $\mathcal{H}_\Lambda$ by
$\mathcal{H}_\lambda$.

We fill a diagram $\lambda$ by hooks to form a tableau $\Lambda$
in the following way: For the first principal hook we fill the
column with entries $1$, $2$, \dots , $r$ and then fill the row
with entries $r+1$, $r+2$, \dots , $s$. We then fill the column of
the second principal hook with $s+1$, $s+2$, \dots , $t$ and fill
the row with $t+1$, $t+2$, \dots , $x$. Continuing in this way we
form the hook tableau.

{\addtocounter{definition}{1} \bf Example \thedefinition .} On the
left is the hook tableau of the diagram $\lambda = (3,3,2)$, and
on the right the same diagram with the content of each box.

\begin{normalsize}

\begin{center}
\begin{picture}(50,50)
\put(0,30){\framebox(15,15)[r]{ 1 }}
\put(15,30){\framebox(15,15)[r]{ 4 }}
\put(30,30){\framebox(15,15)[r]{ 5 }}
\put(0,15){\framebox(15,15)[r]{ 2 }}
\put(15,15){\framebox(15,15)[r]{ 6 }}
\put(30,15){\framebox(15,15)[r]{ 8 }}
\put(0,0){\framebox(15,15)[r]{ 3 }}
\put(15,0){\framebox(15,15)[r]{ 7 }}
\end{picture}
\qquad \qquad \qquad
\begin{picture}(50,50)
\put(0,30){\framebox(15,15)[r]{ 0 }}
\put(15,30){\framebox(15,15)[r]{ 1 }}
\put(30,30){\framebox(15,15)[r]{ 2 }}
\put(0,15){\framebox(15,15)[r]{ -1 }}
\put(15,15){\framebox(15,15)[r]{ 0 }}
\put(30,15){\framebox(15,15)[r]{ 1 }}
\put(0,0){\framebox(15,15)[r]{ -2 }}
\put(15,0){\framebox(15,15)[r]{ -1 }}
\end{picture}
\end{center}
\end{normalsize}
Therefore the sequence $(c_1, c_2, \dots , c_8)$ is given by $(0,
-1, -2, 1, 2, 0 , -1, 1)$. {\nolinebreak \hfill \rule{2mm}{2mm}

Consider (\ref{bigf}) as a rational function of the variables
$z_1, \dots , z_n$ with values in $\mathbb{C}S_n$. The factor
$f_{pq}(z_p + c_p, z_q + c_q)$ has a pole at $z_p = z_q$ if and
only if the numbers $p$ and $q$ stand on the same diagonal of the
tableau $\Lambda$. We then call the pair $(p, q)$ a
\emph{singularity}. And we call the corresponding term $f_{pq}(z_p
+ c_p, z_q + c_q)$ a \emph{singularity term}, or
simply a singularity. \\
Let $p$ and $q$ be in the same hook of $\Lambda$. If $p$ and $q$
are next to one another in the column of the hook then, on
$\mathcal{H}_\lambda$, $f_{pq}(z_p + c_p, z_q + c_q) = 1 - (pq)$.
And so $\frac{1}{2} f_{pq}(z_p + c_p, z_q + c_q)$ is an
idempotent. Similarly, if $p$ and $q$ are next to one another in
the same row of the hook then $f_{pq}(z_p + c_p, z_q + c_q) = 1 +
(pq)$, and
$\frac{1}{2} f_{pq}(z_p + c_p, z_q + c_q)$ will be an idempotent.\\

Also, for distinct $p$, $q$, we have
\begin{equation}\label{inverses} f_{pq}(u,v)f_{qp}(v,u) = 1 - \frac{1}{(u-v)^2}.
\end{equation}
Therefore, if the contents $c_p$ and $c_q$ differ by a number
greater than one, then the factor $f_{pq}(z_p + c_p, z_q + c_q)$
is invertible in $\mathbb{C}S_n$ when $z_p = z_q$ with inverse
$\frac{(c_p - c_q)^2}{(c_p - c_q)^2 - 1} f_{qp}(c_q, c_p)$.

The presence of singularity terms in the product $F_\lambda(z_1,
\dots , z_n)$ mean this product may or may not be regular on the
vector subspace of $\mathcal{H}_\lambda$ consisting of all tuples
$(z_1, \dots , z_n)$ such that $z_1 = z_2 = \dots = z_n$. Using
the following lemma, we will be able to show that $F_\lambda(z_1,
\dots , z_n)$ is in fact regular on this subspace.

\begin{lemma}\label{regular} The restriction of $f_{pq}(u, v)f_{pr}(u, w)f_{qr}(v,
w)$ to the set $(u,v,w)$ such that $u = v \pm 1$ is regular at $u
= w$. \end{lemma}
\begin{proof} Under the condition $u = v \pm 1$, the product $f_{pq}(u, v)f_{pr}(u, w)f_{qr}(v,
w)$ can be written as \[ (1 \mp (pq)) \cdot \left(1 - \frac{(pr) +
(qr)}{v-w}\right). \] And this rational function of $v$, $w$ is
clearly regular at $w = v \pm 1$. \end{proof}

Notice that the three term product, or \emph{triple}, in the
statement of the lemma can be written in reverse order due to
(\ref{triple}). In particular, if the middle term is a singularity
and the other two terms are an appropriate idempotent and
\emph{triple term}, then the triple is regular at $z_1 = z_2 =
\dots = z_n$. We may now prove the first statement of Theorem
\ref{fulltheorem}.

\begin{proposition}\label{jimtheorem1} The restriction of the rational function
$F_\lambda (z_1, \dots , z_n)$ to the subspace
$\mathcal{H}_\lambda$ is regular at $z_1 = z_2 = \dots = z_n$.
\end{proposition}

\begin{proof}
We will prove the statement by reordering the factors of the
product $F_\lambda (z_1, \dots , z_n)$, using relations
(\ref{triple}) and (\ref{commute}), in such a way that each
singularity is part of a triple which is regular at $z_1 = z_2 =
\dots = z_n$, and hence the whole of $F_\lambda (z_1, \dots ,
z_n)$ will be manifestly regular.

Define $g_{pq}$ to be the following; \[ g_{pq} = \left\{
\begin{array}{ccc}
  f_{pq} (z_p + c_p, z_q + c_q) & \textrm{if} & p<q \\
  1 & \textrm{if} & p \geqslant q
\end{array} \right. \]

Now, let us divide the diagram $\lambda$ into two parts,
consisting of those boxes with positive contents and those with
non-positive contents as in Figure \ref{steps}. Consider the
entries of the $i^{\mbox{th}}$ column of the hook tableau
$\Lambda$ of shape $\lambda$ that lie below the steps. If $u_1,
u_2, \dots , u_k$ are the entries of the $i^{\mbox{th}}$ column
below the steps, we define
\begin{equation}\label{cproduct} C_i = \prod_{q=1}^n g_{u_1 , q}
g_{u_2 , q} \dots g_{u_k , q}.
\end{equation} Now consider the entries of the $i^{\mbox{th}}$ row of $\Lambda$ that lie
above the steps. If $v_1, v_2, \dots , v_l$ are the entries of the
$i^{\mbox{th}}$ row above the steps, we define
\begin{equation}\label{rproduct} R_i = \prod_{q=1}^n g_{v_1 , q}
g_{v_2 , q} \dots g_{v_l , q}.
\end{equation}

Our choice of the hook tableau was such that the following is
true; if $d$ is the number of principal hooks of $\lambda$ then by
relations (\ref{triple}) and (\ref{commute}) we may reorder the
factors of $F_\lambda (z_1, \dots , z_n)$ such that \[ F_\lambda
(z_1, \dots , z_n) = \prod_{i=1}^d C_i R_i . \]

Now, each singularity $(p,q)$ has its corresponding term $f_{pq}
(z_p + c_p, z_q + c_q)$ contain in some product $C_i$ or $R_i$.
This singularity term will be on the immediate left of the
 term $f_{p+1,q} (z_{p+1} +
c_{p+1}, z_q + c_q)$. Also, this ordering has been chosen such
that the product of factors to the left of any such singularity in
$C_i$ or $R_i$ is divisible on the right by $f_{p, p+1} (z_p +
c_p, z_{p+1} + c_{p+1})$.\\
Therefore we can replace the pair $f_{pq}(z_p + c_p, z_q + c_q)
f_{p+1,q}(z_p + c_p, z_q + c_q)$ in the product by the triple
$\frac{1}{2} f_{p, p+1} (z_p + c_p, z_{p+1} + c_{p+1}) f_{pq}(z_p
+c_p, z_q + c_q) f_{p+1, q} (z_{p+1} +c_{p+1}, z_q + c_q)$, where
$\frac{1}{2} f_{p, p+1} (z_p + c_p, z_{p+1} + c_{p+1})$ is an
idempotent. Divisibility on the right by $f_{p, p+1}(z_p + c_p,
z_{p+1} + c_{p+1})$ means the addition of the idempotent has no
effect on the value of the product $C_i$ or
$R_i$.\\
By Lemma \ref{regular}, the above triples are regular at $z_1 =
z_2 = \dots = z_n$, and therefore, so too are the products $C_i$
and $R_i$, for all $1 \leqslant i \leqslant d$. Moreover, this
means $F_\lambda (z_1, \dots , z_n)$ is regular at $z_1 = z_2 =
\dots = z_n$.
\end{proof}

{\addtocounter{definition}{1} \bf Example \thedefinition .} As an
example consider the hook tableau of the Young diagram $\lambda =
(3,3,2)$, as shown in Example 2.1.

In the original lexicographic ordering the product $F_\lambda(z_1,
\dots , z_n)$ is written as follows, for simplicity we will write
$f_{pq}$ in place of the term $f_{pq}(z_p + c_p, z_q + c_q)$.
\[
\begin{array}{rl}
   F_\lambda (z_1, \dots , z_n) = & f_{12}f_{13}f_{14}f_{15}\textbf{$f_{16}$}f_{17}f_{18}f_{23}f_{24}f_{25}f_{26}\textbf{$f_{27}$}f_{28}f_{34}f_{35}f_{36} f_{37}f_{38}\\&
   f_{45}f_{46}f_{47}\textbf{$f_{48}$}f_{56}f_{57}f_{58}f_{67}f_{68}f_{78}\\
\end{array}
\]
We may now reorder this product into the form below using
relations (\ref{triple}) and (\ref{commute}) as described in the
above proposition. The terms bracketed are the singularity terms
with their appropriate triple terms.
\[
\begin{array}{rl}
   F_\lambda (z_1, \dots , z_n) = &
   f_{12}f_{13}f_{23}f_{14}f_{24}f_{34}f_{15}f_{25}f_{35}(f_{16}f_{26})f_{36}f_{17}(f_{27}f_{37})f_{18}f_{28}f_{38}\\&
   \cdot f_{45}f_{46}f_{56}f_{47}f_{57}(f_{48}f_{58}) \cdot f_{67}f_{68}f_{78}\\
\end{array}
\]
And so, for each singularity $f_{pq}$, we can replace
$f_{pq}f_{p+1,q}$ in the product by the triple $\frac{1}{2} f_{p,
p+1}f_{pq}f_{p+1, q}$, where $\frac{1}{2} f_{p, p+1}$ is an
idempotent, without changing the value of $F_\lambda (z_1, \dots ,
z_n)$.
\[
\begin{array}{rl}
   F_\lambda (z_1, \dots , z_n) = & f_{12}f_{13}f_{23}f_{14}f_{24}f_{34}f_{15}f_{25}f_{35}(\frac{1}{2}f_{12}f_{16}f_{26})f_{36}f_{17}(\frac{1}{2}f_{23}f_{27}f_{37})f_{18}f_{28}f_{38}\\
   & \cdot f_{45}f_{46}f_{56}f_{47}f_{57}(\frac{1}{2}f_{45}f_{48}f_{58}) \cdot f_{67}f_{68}f_{78}\\
\end{array}
\]
And since each of these triples are regular at $z_1 = z_2 = \dots
= z_n$ then so too is the whole of $F_\lambda (z_1, \dots , z_n)$.
{\nolinebreak \hfill \rule{2mm}{2mm}

\quad

Therefore, due to the above proposition an element $F_\lambda \in
\mathbb{C}S_n$ can now be defined as the value of $F_\lambda (z_1,
\dots , z_n)$ at $z_1 = z_2 = \dots = z_n$. We proceed by showing
this $F_\lambda$ is indeed the diagonal matrix element. To this
end we will need the following propositions.

Note that for $n=1$, we have $F_\lambda =1$. For any $n \geqslant
1$, we have the following fact.

\begin{proposition} The coefficient of $F_\lambda \in
\mathbb{C}S_n$ at the unit element of $S_n$ is 1.
\end{proposition}
\begin{proof} For each $r = 1, \dots , n-1$ let $s_r \in S_n$ be
the adjacent transposition $(r, r+1)$. Let $w_0 \in S_n$ be the
element of maximal length. Multiply the ordered product
(\ref{bigf}) by the element $w_0$ on the right. Using the reduced
decomposition
\begin{equation}\label{longestelement} w_0 = \prod_{(p,q)}^\rightarrow s_{q-p}
\end{equation} we get the product \[ \prod_{(p,q)}^\rightarrow
\left(s_{q-p} - \frac{1}{z_p - z_q + c_p - c_q}\right) . \] For
the derivation of this formula see \cite[(2.4)]{5aa}. This formula
expands as a sum of the elements $s \in S_n$ with coefficients
from the field of rational functions of $z_1, \dots , z_n$ valued
in $\mathbb{C}$. Since the decomposition (\ref{longestelement}) is
reduced, the coefficient at $w_0$ is 1. By the definition of
$F_\lambda$, this implies that the coefficient of $F_\lambda w_0
\in \mathbb{C}S_n$ at $w_0$ is also 1.
\end{proof}

In particular this shows that $F_\lambda \neq 0$ for any nonempty
diagram $\lambda$. Let us now denote by $\varphi$ the involutive
antiautomorphism of the group ring $\mathbb{C}S_n$ defined by
$\varphi (g) = g^{-1}$ for every $g \in S_n$.

\begin{proposition}\label{varphi} The element $F_\lambda \in \mathbb{C}S_n$ is
$\varphi$-invariant. \end{proposition}
\begin{proof} Any element of the group ring $\mathbb{C}S_n$ of
the form $f_{pq}(u,v)$ is $\varphi$-invariant. Applying the
antiautomorphism $\varphi$ to an element of the form (\ref{bigf})
just reverses the ordering of the factors corresponding to the
pairs (p,q). However, the initial ordering can then be restored
using relations (\ref{triple}) and (\ref{commute}), for more
detail see \cite{5aa}. Therefore, any value of the function
$F_\lambda (z_1, \dots , z_n)$ is $\varphi$-invariant. Therefore,
so too is the element $F_\lambda$.
\end{proof}

\begin{proposition}\label{stripping} Let $x$ be last entry in the $k^{\mbox{th}}$ row
of the hook tableau of shape $\lambda$. Denote by $\sigma_x$ the
embedding $\mathbb{C}S_{n-x} \rightarrow \mathbb{C}S_n$ defined by
$ \sigma_x : (pq) \mapsto (p+x, q+x)$ for all distinct $p, q = 1,
\dots , n-x$.\\
 If $\lambda= (\alpha_1, \alpha_2, \dots, \alpha_d
| \beta_1, \beta_2, \dots , \beta_d)$ and $\mu = (\alpha_{k+1},
\alpha_{k+2}, \dots , \alpha_d | \beta_{k+1}, \beta_{k+2}, \dots ,
\beta_d)$, then $F_\lambda = P \cdot \sigma_x(F_\mu)$, for some
element $P \in \mathbb{C}S_n$.
\end{proposition}
\begin{proof}
Here the shape $\mu$ is obtained by removing the first $k$
principal hooks of $\lambda$. By the ordering described in
Proposition \ref{jimtheorem1},
\[ F_\lambda(z_1, \dots ,
z_n) = \prod_{i=1}^k C_iR_i \cdot \sigma_x(F_\mu(z_{x+1}, \dots,
z_{n})),\]
where $C_i$ and $R_i$ are defined by (\ref{cproduct}) and (\ref{rproduct}).\\
Since all products $C_i$ and $R_i$ are regular at $z_1 = z_2 =
\dots = z_n$, Proposition \ref{jimtheorem1} then gives us the
required statement. \end{proof}

In any given ordering of $F_\lambda(z_1, \dots , z_n)$, we want a
singularity term to be placed next to an appropriate triple term
such that we may then form a regular triple. In that case we will
say these two terms are 'tied'.\\
We now complete the proof of Theorem \ref{fulltheorem}. If $u$,
$v$ stand next to each other in the same row, or same column, of
$\Lambda$ the following results show that $F_\lambda$ is divisible
on the left (and on the right) by $1 - (uv)$, or $1+(uv)$
respectively, or divisible by these terms preceded (followed) by
some invertible elements of $\mathbb{C}S_n$. Hence $F_\lambda$ is
the diagonal matrix element. \\
However, proving the divisibilities described requires some pairs
to be 'untied', in which case we must form a new ordering. This is
the content of the following proofs. Some explicit examples will
then given in Example 2.10 below.

\begin{proposition}\label{jimtheorem2} Suppose the numbers $u < v$ stand next to each
other in the same column of the hook tableau $\Lambda$ of shape
$\lambda$. First, let $s$ be the last entry in the row containing
$u$. If $c_v < 0$ then the element $F_\lambda \in \mathbb{C}S_n$
is divisible on the left and on the right by $f_{u,v}(c_u, c_v) =
1 - (uv)$. If $c_v \geqslant 0$ then the element $F_\lambda \in
\mathbb{C}S_n$ is divisible on the left by the product
\[ \prod_{i = u, \dots, s}^\leftarrow \left( \prod_{j= s+1, \dots,
v}^\rightarrow f_{ij}(c_i, c_j) \right) \]
\end{proposition}

\begin{proof}
By Proposition \ref{varphi}, the divisibility of $F_\lambda$ by
the element $1- (uv)$ on the left is equivalent to the
divisibility by the same element on the right. Let us prove
divisibility on the left.

By Proposition \ref{stripping}, if $\sigma_x(F_\mu)$ is divisible
on the right by $f_{uv}(c_u, c_v)$, or $f_{uv}(c_u, c_v)$ followed
by some invertible terms, then so too will $F_\lambda$. If
$\sigma_x(F_\mu)$ is divisible on the left by $f_{uv}(c_u, c_v)$,
or $f_{uv}(c_u, c_v)$ preceded by some invertible terms, then, by
using Proposition \ref{varphi} twice, first on $\sigma_x(F_\mu)$
then on the product $F_\lambda = P \cdot \sigma_x(F_\mu)$, so too
will $F_\lambda$. Hence we only need to prove the statement for
$(u,v)$ such that $u$ is in the first row or first column of
$\Lambda$.

Let $r$ be the last entry in the first column of $\Lambda$, $s$
the last entry in the first row of $\Lambda$, and $t$ the last
entry in the second column of $\Lambda$, as shown in Figure
\ref{jimtheorem2fig}.

\begin{figure}[h] \label{jimtheorem2fig}

\begin{normalsize}
\begin{center}
\begin{picture}(275,200)
\begin{small}
\put(0,175){\framebox(25,25)[c]{ 1 }}
\put(25,175){\framebox(25,25)[c]{ $r + 1$ }}
\put(50,175){\framebox(25,25)[c]{ $r+2$ }}
\put(0,150){\framebox(25,25)[c]{ 2 }}
\put(25,150){\framebox(25,25)[c]{ $s+1$ }}
\put(50,150){\framebox(25,25)[c]{ $t+1$ }}

\put(0,25){\line(0,1){125}} \put(25,25){\line(0,1){125}}
\put(0,0){\framebox(25,25)[c]{ $r$ }}

\put(75,200){\line(1,0){50}} \put(75,175){\line(1,0){50}}
\put(125,175){\framebox(25,25)[c]{ $u$ }}
\put(150,200){\line(1,0){100}} \put(150,175){\line(1,0){100}}
\put(250,175){\framebox(25,25)[c]{ $s$ }}

\put(50,100){\line(0,1){50}} \put(25,75){\framebox(25,25)[c]{ $t$
}}

\put(75,150){\line(1,0){50}} \put(125,150){\framebox(25,25)[c]{
$v$ }} \put(150,150){\line(1,0){50}} \put(200,150){\line(0,1){25}}

\linethickness{1.5pt} \put(0,200){\line(1,0){25}}
\put(25,175){\line(1,0){25}} \put(50,150){\line(1,0){25}}
\put(25,175){\line(0,1){25}} \put(50,150){\line(0,1){25}}

\put(125,150){\line(0,1){50}}\put(150,150){\line(0,1){50}}
\put(125,200){\line(1,0){25}}\put(125,150){\line(1,0){25}}

\put(95,185){$\dots$} \put(170,185){$\dots$}\put(220,185){$\dots$}
\put(95,160){$\dots$} \put(170,160){$\dots$}

\put(10,120){$\vdots$} \put(35,
120){$\vdots$}\put(10,80){$\vdots$} \put(10,40){$\vdots$}

\end{small}
\end{picture}
\end{center}
\end{normalsize}
 \caption{The first two principal hooks of the hook tableau $\Lambda$}
 \end{figure}

We now continue this proof by considering three cases and showing
the appropriate divisibility in each.

\emph{(i)} \quad  If $c_v < 0$ (i.e. $u$ and $v$ are in the first
column of $\Lambda$) then $v = u+1$ and $F_\lambda(z_1, \dots ,
z_n)$ can be written as $F_\lambda(z_1, \dots, z_n) = f_{uv}(z_u
+c_u, z_v
+c_v) \cdot F$. \\
Starting with $F_\lambda(z_1, \dots , z_n)$ written in the
ordering described in Proposition \ref{jimtheorem1} and simply
moving the term $f_{u,v}(z_u + c_u, z_v + c_v)$ to the left
results in all the singularity terms in the product $F$ remaining
tied to the same triple terms as originally described in that
ordering. Therefore we may still form regular triples for
each singularity in $F$, and hence $F$ is regular at $z_1 = z_2 = \dots = z_n$.\\
So by considering this expression for $F_\lambda(z_1, \dots ,
z_n)$ at $z_1 = z_2 = \dots = z_n$ we see that $F_\lambda$ will be
divisible on the left/right by $f_{uv}(c_u, c_v) = (1 - (uv))$.

\emph{(ii)} \quad  If $c_v = 0$ then $v=s+1$, and $F_\lambda(z_1,
\dots , z_n)$ can be written as \[ F_\lambda(z_1, \dots , z_n) =
\prod_{i = u, \dots , s}^\leftarrow f_{i, s+1}(z_i + c_i, z_{s+1}
+ c_{s+1}) \cdot F' . \] Again, starting with the ordering
described in Proposition \ref{jimtheorem1}, this results in  all
the singularity terms in the product $F'$ remaining tied to the
same triple terms as originally described in that ordering. Hence
$F'$ is regular at $z_1 = z_2 = \dots = z_n$. And so $F_\lambda$
is divisible on the left by
\[ \prod_{i = u, \dots , s}^\leftarrow f_{i, s+1}(c_i, c_{s+1}).
\]

\emph{(iii)} \quad  If $c_v > 0$ (i.e. $v$ is above the steps)
then $f_{uv}(z_u + c_u, z_v + c_v)$ is tied to the singularity
$f_{u-1, v}(z_{u-1} + c_{u-1}, z_v + c_v)$ as a triple term. To
show divisibility by $f_{uv}(z_u + c_u, z_v + c_v)$ in this case
we need an alternative expression for $F_\lambda(z_1, \dots ,
z_n)$ that is regular when $z_1 = z_2 = \dots = z_n$. Define a
permutation $\tau$ as follows,

\[ \tau = \prod_{i = u, \dots, s}^\rightarrow \left( \prod_{j= s+1, \dots,
v}^\leftarrow (i j) \right)
\phantom{XXXXXXXXXXXXXXXXXXXXXXXXXXXXXXX}
\]
\[ = \left(
\small \begin{array}{ccccccccccccccccc}
                   1 & 2 & \dots & u-1 & u & u+1 & \dots &  & \dots & \dots  &  & \dots & v-1 & v & v+1 & \dots & n \\
                   1 & 2 & \dots & u-1 & s+1 & s+2 & \dots & v-1 & v & u & u+1 & \dots & s-1 & s & v+1 & \dots & n \\
                 \end{array} \right) \]

From the definition of $C_1$ in (\ref{cproduct}) we now define
$C'_1 = \tau C_1$, where $\tau$ acts on the indices of the product
$C_1$ such that
\[ \tau C_1 = \prod_{j=1}^n \left( \prod_{i=1}^r g_{i, \tau
\cdot j} \right). \]

For the rest of this proof we will simply write $f_{ij}$ instead
of $f_{ij}(z_i + c_i, z_j + c_j)$. Define $R'_1$ as,
\begin{eqnarray*}
R'_1 & = & \prod_{i= r+2, \dots, u-1}^\leftarrow \left(
\prod_{j=s+1, \dots ,v}^\rightarrow f_{ij} \right) \cdot \prod_{i=
s+1, \dots, t-1}^\rightarrow \left( \prod_{j=i+1, \dots
,t}^\rightarrow f_{ij} \right) \cdot \left( \prod_{j=s+1, \dots
,t}^\leftarrow f_{r+1, j} \right) \left( \prod_{j=t+1, \dots
,v}^\rightarrow f_{r+1, j} \right) \\
   && \times \prod_{i= r+1,
\dots, s-1}^\rightarrow \left( \prod_{j=i+1, \dots ,s}^\rightarrow
f_{ij} \right) \cdot \prod_{j= v+1, \dots, n}^\rightarrow
\left( \prod_{i=r+1, \dots ,s}^\rightarrow f_{ij} \right). \\
\end{eqnarray*}

Finally, define $C'_2$ as, \[ C'_2 = \prod_{j= t+1, \dots,
n}^\rightarrow \left( \prod_{i=s+1, \dots ,t}^\rightarrow f_{ij}
\right). \]

Then, \[ F_\lambda(z_1, \dots , z_n) = \prod_{i = u, \dots ,
s}^\leftarrow \left( \prod_{j=s+1, \dots, v}^\rightarrow f_{ij}
\right) \cdot C'_1 R'_1 C'_2 R_2 \cdot \prod_{i=3}^d C_iR_i, \]
where $d$ is the number of principal hooks of $\lambda$.

The product $C'_1 R'_1 C'_2 R_2 \cdot \prod C_iR_i$ is regular at
$z_1 = z_2 = \dots z_n$ since, as before, for any singularity
$(p,q)$ the terms $f_{pq}f_{p+1 ,q}$ can be replaced by the triple
$\frac{1}{2}f_{p, p+1}f_{pq}f_{p+1, q}$ -- except in the
expression $R'_1$ where the terms $f_{pl}f_{pq}$ are replaced by
$\frac{1}{2} f_{pl}f_{pq}f_{lq}$, where $l$ is the entry to the
immediate left of $q$. Note that $l = q-1$ when $c_q > 1$ and $l =
s+1$ when $c_q = 1$. \\
And so by letting $z_1 = z_2 = \dots =z_n$ we see that $F_\lambda$
is divisible on the left by \[ \prod_{i = u, \dots , s}^\leftarrow
\left( \prod_{j=s+1, \dots, v}^\rightarrow f_{ij}(c_i, c_j)
\right).
\]
\end{proof}

\begin{proposition}\label{jimtheorem3} Suppose the numbers $u < v$ stand next to each
other in the same row of the hook tableau $\Lambda$ of shape
$\lambda$. Let $r$ be the last entry in the column containing $u$.
If $c_u > 0$ then the element $F_\lambda \in \mathbb{C}S_n$ is
divisible on the left and on the right by $f_{u,v}(c_u, c_v) = 1 +
(uv)$. If $c_u \leqslant 0$ then the element $F_\lambda \in
\mathbb{C}S_n$ is divisible on the left by the product
\[ \prod_{i = u, \dots, r}^\leftarrow \left( \prod_{j= r+1, \dots,
v}^\rightarrow f_{ij}(c_i, c_j) \right) \]
\end{proposition}

We omit the proof of this proposition as it is very similar to
that of Proposition \ref{jimtheorem2}. 

This completes our proof of Theorem \ref{fulltheorem}. Let us now
consider an example that allows us to see how the product
$F_\lambda(z_1, \dots, z_n)$ is broken down in the proof of
Proposition \ref{jimtheorem2}.

{\addtocounter{definition}{1} \bf Example \thedefinition .} We
again consider the hook tableau of the Young diagram $\lambda =
(3,3,2)$, as shown in Example 2.1.

We begin with the product $F_\lambda(z_1, \dots, z_n)$ in the
ordering described in Proposition \ref{jimtheorem1}. For
simplicity we again write $f_{pq}$ in place of the term
$f_{pq}(z_p + c_p, z_q + c_q)$. We have also marked out the
singularities in this expansion along with their triple terms, but
no idempotents have yet been added which would form regular
triples.

\begin{equation}\label{example1}
\begin{array}{rl}
   F_\lambda (z_1, \dots , z_n) = &
   f_{12}f_{13}f_{23}f_{14}f_{24}f_{34}f_{15}f_{25}f_{35}(f_{16}f_{26})f_{36}f_{17}(f_{27}f_{37})f_{18}f_{28}f_{38}\\&
   \cdot f_{45}f_{46}f_{56}f_{47}f_{57}(f_{48}f_{58}) \cdot f_{67}f_{68}f_{78}\\
\end{array}
\end{equation}

Let $u=4$ and $v=6$ in Proposition \ref{jimtheorem2}. Then by that
proposition we may arrange the above product as follows. Notice
since $c_v = 0$ all singularity-triple term pairs remain the same.

\[
\begin{array}{rl}
   F_\lambda (z_1, \dots , z_n) = &
   f_{56}f_{46} \cdot f_{12}f_{13}f_{23}(f_{16}f_{26})f_{36}f_{14}f_{24}f_{34}f_{15}f_{25}f_{35}f_{17}(f_{27}f_{37})\\&
   f_{18}f_{28}f_{38} \cdot f_{45}f_{47}f_{57}(f_{48}f_{58}) \cdot f_{67}f_{68}f_{78}\\
\end{array}
\]
We may now add the appropriate idempotents so that all
singularities remain in regular triples. And so by considering the
product at $z_1 = z_2 = \dots = z_n$ we have that $F_\lambda$ is
divisible on the left by $(1 - (46))$, preceded only by invertible
terms, as desired.

Now let $u=5$ and $v=8$. As described by Proposition
\ref{jimtheorem2} (iii) we may arrange (\ref{example1}) in the
following way. Singularities in $R_1$ have been marked out with
their alternative triple terms, while all other singularity-triple
term pairs remain the same.
\[
\begin{array}{rl}
   F_\lambda (z_1, \dots , z_n) = &
   f_{56}f_{57}f_{58} \cdot f_{12}f_{13}f_{23}f_{14}f_{24}f_{34}(f_{16}f_{26})f_{36}f_{17}(f_{27}f_{37})f_{18}f_{28}f_{38}\\&
   f_{15}f_{25}f_{35} \cdot f_{67}f_{47}(f_{46}f_{48})f_{45} \cdot f_{68}f_{78}\\
\end{array}
\]
In moving $f_{58}$ to the left it is untied from the singularity
$f_{48}$. So we must form new triples which are regular at $z_1 =
z_2 = \dots = z_n$. Therefore, by considering the product at $z_1
= z_2 = \dots = z_n$, we have that $F_\lambda$ is divisible on the
left by $(1 - (58))$, again preceded only by invertible terms, as
desired. {\nolinebreak \hfill \rule{2mm}{2mm}

\end{document}

%% file: headerejc.tex
\newenvironment{proof}{\begin{trivlist} \item[]
{\bf Proof.}}{\nolinebreak
\hfill \rule{2mm}{2mm} \end{trivlist}}

\newcounter{ctr}

\newcounter{ctr1}

\newcounter{ctr2}

\newcounter{ctr3}

\newtheorem{definition}{Definition}[section]    
\newtheorem{theorem}[definition]{Theorem}
\newtheorem{nottheorem}[definition]{}
\newtheorem{lemma}[definition]{Lemma}

\newtheorem{proposition}[definition]{Proposition}

\newcommand{\NN}{\hbox{{\sf I}\kern-.15em\hbox{\sf N}}}
\newcommand{\ZZ}{\hbox{{\sf Z}\kern-.77em\hbox{\sf Z}\kern.15em}}
\newcommand{\QQ}{\hbox{{\sf I}\kern-.4em\hbox{\sf Q}}}
\newcommand{\RR}{\hbox{{\sf I}\kern-.15em\hbox{\sf R}}}
\newcommand{\CC}{\hbox{{\sf I}\kern-.4em\hbox{\sf C}}}

\newcommand{\ud}{\, {\rm d} \kern-.015em }

\newcommand{\mod}[1]{\left| \kern.05em #1 \kern.05em \right|}
\newcommand{\norm}[1]{\left\| \kern.05em #1 \kern.05em \right\|}
\newcommand{\inner}[1]{\left\langle \kern.05em #1 \kern.05em \right\rangle }

\newcommand{\pick}[2]{\renewcommand{\arraystretch}{0.6}
\left( \kern-.4em \begin{array}{c} #1 \\ #2 \end{array} \kern-.4em \right) }

